\documentclass[a4paper, 12pt]{article}
\usepackage[margin=1in]{geometry}

\usepackage{amsmath}
\usepackage{amsfonts}
\usepackage{amssymb}
\usepackage{amsthm}
\usepackage{bbm}
\usepackage{mathrsfs}
\usepackage{mathtools}
\usepackage{comment}
\usepackage{dsfont}
\usepackage[dvipsnames]{xcolor}
\usepackage[normalem]{ulem}
\usepackage{hyperref}
\usepackage{listings}
\usepackage{float}

\usepackage{bm}
\usepackage{bbm}
\usepackage{natbib}
\definecolor{darkblue}{rgb}{0.0,0.0,0.65}
\definecolor{darkred}{rgb}{0.65,0.0,0.0}
\hypersetup{
  colorlinks = true,
  citecolor  = darkred,
  linkcolor  = darkblue,
  filecolor  = darkblue,
  urlcolor   = darkblue,
}
\AtBeginDocument{
    \hypersetup{colorlinks,urlcolor=black}
}

\usepackage{thmtools, thm-restate}

\newcommand{\pref}[1]{(\ref{#1})}

\let\temp\phi
\let\phi\varphi
\let\varphi\temp

\let\temp\epsilon
\let\epsilon\varepsilon
\let\varepsilon\temp

\DeclareMathOperator{\Id}{Id}

\DeclareMathOperator{\Div}{div}
\DeclareMathOperator{\clos}{clos}
\DeclareMathOperator{\WFR}{WFR}
\DeclareMathOperator{\D}{D}

\DeclareMathOperator*{\argmin}{arg\,min}

\DeclareMathOperator{\textKL}{KL}

\DeclarePairedDelimiterX{\infdivx}[2]{(}{)}{%
  #1\;\delimsize\|\;#2%
}
\newcommand{\KL}{\textKL\infdivx}

\newtheorem{theorem}{Theorem}

\newtheorem{prop}{Proposition}

\newtheorem*{remark}{Remark}

\begin{document}
\title{Wasserstein-Fisher-Rao Splines}
\date{\today}

\author{
  {\normalsize Julien Clancy}\\
  {\small \tt clancy.julien@gmail.com}
  \and
  {\normalsize Felipe Suárez} \\
  {\small \tt felipesc@mit.edu}
  \and {\normalsize Massachusetts Institute of Technology}\thanks{J.C. was supported by NSF GRFP, and F.S. was suppported by the MathWorks Fellowship and the NSF grant IIS-1838071.}
}

% \author[1]{Julien Clancy \thanks{ Email: {\tt clancy.julien@gmail.com}. This work by NSF GRFP.}}
% \author[2]{Felipe Suarez \thanks{Email: {\tt felipesc@mit.edu}. This work was supported by the MathWorks Fellowship and the NSF grant IIS-1838071.}}
% \affil[2]{Department of Mathematics, MIT}

\maketitle
\abstract{We study interpolating splines on the Wasserstein-Fisher-Rao (WFR) space of measures with differing total masses. To achieve this, we derive the covariant derivative and the curvature of an absolutely continuous curve in the WFR space. We prove that this geometric notion of curvature is equivalent to a Lagrangian notion of curvature in terms of particles on the cone. Finally, we propose a practical algorithm for computing splines extending the work of \cite{chewi2020fast}.}

\section{Introduction}

Let $\mu_1,\mu_2,\dots,\mu_n$ be $n$ positive measures of differing total masses. How to interpolate them? This question is motivated by cellular trajectory reconstruction where $\mu_i$ is a population of cells at time $i$ \citep{schiebinger2019optimal}. Cells move in gene space as they evolve, but also divide and proliferate. While ad-hoc fixes for this issue have been proposed, e.g. via renormalization and using optimal transport (OT) \citep{chewi2020fast}, the conservation of mass property inherent to OT makes it a less suitable tool for this task.

Curve evolution in Wasserstein space is governed by the continuity equation
\begin{equation}\label{w2_ce}
    \partial_t \mu_t + \Div (v_t \mu_t) = 0,
\end{equation}
which can be viewed as a simple restatement of conservation of mass, following the divergence theorem.\footnote{Formally, this equation is to be interpreted weakly in duality with functions $\varphi \in \mathcal{C}^\infty_c(\mathbb{R}^d \times \mathbb{R})$ via the divergence theorem i.e. $\tfrac{d}{dt} \int \phi\,d\mu_t = \int \langle v_t, \nabla \phi \rangle \, d\mu_t$ .} This equation is underdetermined, and the geometry of $W_2$ is induced by selecting for each time $t$ the field $v_t$ with minimal kinetic energy:
\begin{equation*}
    v_t = \argmin_{u_t} \int \lVert u_t \rVert^2 \, d \mu_t = \argmin_{u_t} \, \lVert u_t \rVert_{\mu_t}^2.
\end{equation*}

It can be seen \citep{gigli2012second} that the optimal $v_t$ lies in the closure of $\{\nabla \phi \mid \phi \in C^\infty_c\}$. The celebrated Benamou-Brenier theorem states that the $W_2$ distance, defined by optimal transport, is equal to the least total kinetic energy among all possible paths.
\begin{theorem}\citep{benamou2000computational}
    Let $\mu_0,\mu_1$ be probability measures. Then
    \begin{equation*}
        W_2^2(\mu_0,\mu_1) = \inf_{(\mu_t,v_t)} \int_0^1 \lVert v_t \rVert_{\mu_t}^2 \, dt,
    \end{equation*}
    where the infimum is taken over solutions of the continuity equations with prescribed boundary data $\mu_0$ and $\mu_1$.
\end{theorem}

Thus, in order to define a new metric between measures of arbitrary mass we have to alter the continuity equation, running the procedure in reverse. Turning back to intuition, where the term $\Div (v_t \mu_t)$ represents mass translation, we add another term representing growth or decay:
\begin{equation}\label{eq:wfr_ce}
    \partial_t \mu_t + \Div(v_t \mu_t) = 4 \alpha_t \mu_t.
\end{equation}
We call this the {\it nonconservative continuity equation}.\footnote{The factor of $4$ is for notational convenience in accordance with the literature.} The term in the right-hand side allows for a \textit{relative} growth or decay in mass. While measures $(\mu_t)$ evolving according to \eqref{eq:wfr_ce} may have varying mass, they are granted to stay positive.

From here, we measure the magnitude of a pair $(v_t, \alpha_t)$ by
\begin{equation*}
    \lVert (v_t,\alpha_t) \rVert_{\mu_t}^2 = \int \left(| v_t |^2 + 4 \alpha_t^2 \right) \,d \mu_t,
\end{equation*}
where $\| \cdot \|$ stands for norm in WFR space and $|\cdot |$ for norm of vectors in Euclidean space.
The Wasserstein-Fisher-Rao distance is then defined by
\begin{equation}\label{eq:wfr_variational}
    \WFR(\mu_0,\mu_1)^2 = \inf_{(\mu_t,v_t,\alpha_t)} \int_0^1 \lVert (v_t,\alpha_t) \rVert_{\mu_t}^2 \, dt
\end{equation}
where, again, the minimization is over solutions to \pref{eq:wfr_ce} with $\mu_0$ and $\mu_1$ prescribed. For more detailed discussions of the genesis of the equation \pref{eq:wfr_ce} and the development of the connection with optimal transport see \citealp{kondratyev2016new,liero2018optimal,chizat2017unbalanced}.

It was shown simultaneously in \citet{liero2018optimal,chizat2017unbalanced} that \pref{eq:wfr_variational} defines a metric on $\mathcal{M}_+(\mathbb{R}^d)$, the space of non-negative measures, and that it turns this into a geodesic space. Furthermore, \citet{kondratyev2016new} shows it has the structure of a pseudo-Riemannian manifold analogous to the Riemannian structure on $W_2$ \citep{otto2001geometry}, with inner product given by
\begin{equation*}
    \left \langle (v, \alpha), (w, \beta) \right \rangle_{\mu} = \int \left( \langle v, w \rangle + 4 \alpha \beta \right) \, d \mu,
\end{equation*}
where the tangent space is
\begin{equation*}
    T_\mu (\mathcal{M}_+) = \clos_{L^2(\mu)} \left\{ (\nabla \alpha, \alpha) \mid \alpha \in \mathcal{C}^\infty_c(\mathbb{R}^d) \right\}.
\end{equation*}

\begin{remark}
Notice that this is the same tangent space as for $W_2$, with the difference being that the $\WFR$ Riemannian metric is the full $H^1$ Sobolev norm, while the $W_2$ Riemannian metric includes only the gradient.
\end{remark}

There has been recent interest in understanding the curvature of these spaces \citep{chewi2020gradient}, and specifically in understanding curvature-minimizing interpolating curves that generalize Euclidean splines \citep{benamou2019second,chen2018measure,chewi2020fast}. We aim to replicate some of these advances in the Wasserstein-Fisher-Rao case, which, as mentioned, plays an important role in applications. Specifically, we will 1) characterize the covariant derivative in WFR space and use this to define a notion of intrinsic splines; 2) define splines from a pseudo-Eulerian perspective, as measures on paths, and establish a relationship to intrinsic splines; and 3) define a more practically workable notion of WFR splines, analogous to the transport splines of \citet{chewi2020fast}, and examine its relationship to intrinsic splines.

\section{The Covariant Derivative}

In this section, we derive an expression for the covariant derivative. We first present a new derivation of the covariant derivative in $W_2$ that is simpler than the definition in \citet{gigli2012second}. In turn, this new derivation is extended to WFR space.

The Riemannian metric on $W_2$ is given by $\langle v, w \rangle_\mu = \int \langle v, w \rangle \, d \mu$, and the covariant derivative $\frac{\D}{dt}$ and its associated Levi-Civita connection $\nabla$ must satisfy two properties:
\begin{enumerate}
    \item \textbf{Leibniz rule.} If $\mu_t$ is a curve and $(v_t^{1})$ and $(v_t^{2})$ are two (tangent) vector fields along it, then
    \begin{equation}\label{eq:leibniz}
        \frac{d}{dt} \langle v_t^1, v_t^2 \rangle_{\mu_t} = \left \langle v_t^1, \frac{\mathbf{D}}{dt} v_t^2 \right \rangle_{\mu_t} + \left \langle \frac{\mathbf{D}}{dt} v_t^1, v_t^2 \right \rangle_{\mu_t}.\tag{P1}
    \end{equation}
    \item \textbf{Torsion-freeness.} If $X$ and $Y$ are vector fields, then
    \begin{equation}\label{eq:torsion}
        \nabla_X Y - \nabla_Y X = [X, Y]. \tag{P2}
    \end{equation}
\end{enumerate}

Let $(v_t)$ be the tangent field of the curve $(\mu_t)$ (as noted above, the minimal field is unique and is a gradient\footnote{In general it is merely in the closure of the set of gradients, but if all measures involved are absolutely continuous then it is truly a gradient.}). The product rule yields
\begin{align*}
    \frac{d}{dt} \langle v_t^1, v_t^2 \rangle_{\mu_t} &= \frac{d}{dt} \int \langle v_t^2, v_t^2 \rangle \, d \mu_t\\
    &= \int \left( \langle \partial_t v_t^2, v_t^2 \rangle + \langle v_t^2, \partial_t v_t^2 \rangle\right) \, d \mu_t + \int \langle v_t^2, v_t^2 \rangle \, d (\partial_t \mu_t).
\end{align*}
Because $(\mu_t, v_t)$ solves the continuity equation, the dual definition of $\Div (v_t \mu_t)$ and the divergence theorem give
\begin{equation*}
    \int \langle v_t^1, v_t^2 \rangle \, d (\partial \mu_t) = \int \langle \nabla v_t^1 \cdot v_t, v_t^2 \rangle + \langle \nabla v_t^2 \cdot v_t, v_t^1 \rangle \, d \mu_t.
\end{equation*}
Together with \eqref{eq:leibniz}, it yields 
\begin{equation*}
    \left \langle \frac{\mathbf{D}}{dt} v_t^1, v_t^2 \right \rangle_{\mu_t} + \left \langle v_t^1, \frac{\mathbf{D}}{dt} v_t^2 \right \rangle_{\mu_t} = \langle \partial_t v_t^1 + \nabla v_t^1 \cdot v_t, v_t^2 \rangle_{\mu_t} + \langle v_t^1, \partial_t v_t^2 + \nabla v_t^2 \cdot v_t \rangle_{\mu_t}.
\end{equation*} 

From here, it natural to postulate that
\begin{equation}\label{eq:w2_cov_derivative}
    \frac{\mathbf{D}}{dt}v_t^1 = \mathcal{P}_{\mu_t}\left( \frac{D}{dt} v_t^1 \right) = \mathcal{P}_{\mu_t} \left( \partial_t v_t^1 + \nabla v_t^1 \cdot v_t \right),
\end{equation}
where $\mathcal{P}_{\mu_t}$ is the orthogonal projection onto $T_{\mu_t}(\mathcal{P}_2)$ in $L^2(\mu_t)$. We call $\frac{D}{dt}$ the {\it total derivative}. Note that if $v_t^1 = v_t = \nabla \phi_t$ then $\partial_t v_t + \nabla v_t \cdot v_t = \nabla \left(\partial_t \phi_t + \frac{1}{2}|\nabla \phi_t|^2 \right) \in T_{\mu_t}(\mathcal{P}_2)$, so no projection is necessary, and the total and covariant derivatives coincide; in other words,
\begin{equation}\label{eq:w2_d2dt2}
    \frac{\mathbf{D}^2}{dt^2} \mu_t = \partial_t v_t + \nabla v_t \cdot v_t.
\end{equation}

We now examine the torsion-free property \eqref{eq:torsion}, and we follow Gigli's argument \cite[Section~\S 5.1]{gigli2012second}, which we partially repeat for convenience of reference. It is quite technical to give meaning directly to a smooth vector field on all of $\mathcal{P}_2$, and so to the Levi-Civita connection, but it can be indirectly defined via the covariant derivative, which is all that will be necessary in this work. Let $(\mu^{1}_t)$ and $(\mu_t^{2})$ be two absolutely continuous curves with measures that are absolutely continuous with respect to the Lebesgue measure, such that $\mu^{1}_0 = \mu^{2}_0 = \mu$, and let their velocity fields be $(v^{1}_t)$ and $(v^{2}_t)$. Since $v_0^1$, $v_0^2$ are gradients, they are in $T_\mu(\mathcal{P}_2)$ for every $\mu$, so we may define two new tangent fields along these curves by
\begin{align*}
    u_t^{1} &= v^{2}_0,\\
    u_t^{2} &= v^{1}_0.
\end{align*}

With this definition, it is reasonable to interpret
\begin{equation*}
    \nabla_{u^{1}_0} u^{2}_t \Big|_{t=0} = \frac{\mathbf{D}}{dt} u_t^{2} \Big|_{t=0},
\end{equation*}
with the derivative being taken along $\mu_t^{2}$, and similarly for $\nabla_{u^{2}_0} u^{1}_t$. Now, fix $\phi$ and consider the functional $F\colon \mu \mapsto \int \phi \, d \mu$. By the continuity equation, the derivative of $F$ along $u_0^{2}$ at $\mu$ is
\begin{equation*}
    \frac{d}{dt} F[\mu_t^{1}] \Big|_{t=0} = \int \phi \, d \left( \partial_t \mu_t^{1} \right) \Big|_{t=0} = \int \langle \nabla \phi, u_0^{2} \rangle \, d \mu.
\end{equation*}
Then since the covariant derivative above respects the metric,
\begin{align*}
    u_0^{1}(u_0^{2}(F))[\mu] &= \frac{d}{dt} \langle \nabla \phi, u_0^{2} \rangle_{\mu_t^{2}} \Big|_{t=0}\\
    &= \left \langle \frac{\mathbf{D}}{dt} \nabla \phi, u_0^{2} \right \rangle_{\mu_t^{2}} + \left \langle \nabla \phi, \frac{\mathbf{D}}{dt} u_0^{2} \right \rangle_{\mu_t^{2}} \bigg|_{t=0}\\
    &= \left \langle \nabla^2 \phi \cdot u^{1}_0, u^{2}_0 \right \rangle_\mu + \left \langle \nabla \phi, \nabla_{u^{1}_0} u^{2}_t \right \rangle_\mu
\end{align*}
where we have use the definition of $\frac{\mathbf{D}}{dt}$ to calculate $\frac{\mathbf{D}}{dt} \nabla \phi$ on the third line. Performing the same calculation for $u_0^{2}(u_0^{1}(F))[\mu]$ and subtracting, since $\nabla^2 \phi$ is symmetric the first terms cancel and we get
\begin{equation*}
    u_0^{1}(u^{2}(F))[\mu] - u_0^{2}(u^{1}(F))[\mu] = \left \langle \nabla \phi, \nabla_{u^{1}_0} u^{2}_t - \nabla_{u^{2}_0} u^{1}_t \right \rangle_\mu
\end{equation*}
Since gradients $\nabla \phi$ are dense in $T_\mu(\mathcal{P}_2)$ this means that $\frac{\mathbf{D}}{dt}$ is indeed torsion-free.

Now we repeat the argument in WFR space. 

\begin{restatable}{theorem}{firstthm}\label{thm:covariant_statement}
   Let $\mu_t$ be an absolutely continuous curve in Wasserstein-Fisher-Rao space satisfying the continuity equation with tangent fields $(v_t,\alpha_t)$. The covariant derivative is given by
   \begin{equation}\label{eq:wfr_d2dt2}
    \frac{\mathbf{D}^2}{dt^2} \mu_t = \begin{pmatrix}
        \partial_t v_t + \nabla v_t \cdot v_t + 4 \alpha_t v_t\\
        \partial_t \alpha_t + \frac{1}{2}|\nabla \alpha_t|^2 + 2 \alpha_t^2
    \end{pmatrix}.
\end{equation}
\end{restatable}

\begin{proof}
See appendix \ref{thm:covariant_derivative}.
\end{proof}

Theorem \ref{thm:covariant_statement} gives a geometrical proof of the dynamical characterization of geodesics in duality with \eqref{eq:wfr_variational} (\citealp[Theorem~1.1.15]{chizat2017unbalanced}; \citealp[Theorem~8.12]{liero2018optimal}). Indeed, an absolutely continuous geodesic $(\mu_t)$ in WFR space with tangent $\dot{\mu}_t=(\nabla \alpha, \alpha)$  satisfies the Hamilton-Jacobi equation almost surely
\begin{equation*}
    \partial_t \alpha_t + \frac{1}{2}|\nabla \alpha_t|^2 + 2 \alpha_t^2 = 0,
\end{equation*}
or, equivalently, the curve $(\mu_t)$ is autoparallel:
\begin{equation}
    \nabla_{\dot{\mu}_t} \dot{\mu}_t = \frac{\mathbf{D}^2}{dt^2} \mu_t = 0.
\end{equation}

\section{E-Splines and P-splines}
In this section, we define notions of $E$- and $P$-splines. As in the $W_2$ case, we can define an intrisic notion of curvature-minimizing interpolators, which we term {\it $E$-splines}, by
\begin{equation}\label{eq:wfr_espline}
    \inf_{(\mu_t,\mathbf{v}_t)} \int_0^1 \left \lVert \frac{\mathbf{D}}{dt} \mathbf{v}_t \right \rVert_{\mu_t}^2 \, dt \text{ s.t. } \mu_{t_i} = \mu_i
\end{equation}
Though the characterization \pref{eq:wfr_d2dt2} yields an explicit objective function, there is no practical way to optimize \pref{eq:wfr_espline}. Thus, we first define an analogy to the path splines ({\it P-splines}) introduced in \citet{chen2018measure,benamou2019second}. We give a brief introduction to them here (see \citealp{chewi2020fast} for more details).

\subsection{Geodesics and \texorpdfstring{$P$}{P}-splines in \texorpdfstring{$W_2$}{W2}}

In the Wasserstein space $W_2$ over a metric space $(\mathcal{X},d)$ it is known that geodesics can be represented as measures over paths (see \citealp{lisini2007characterization}). Specifically, let $\Omega$ be the set all absolutely continuous paths in $\mathcal{X}$, let $l$ be the length functional on $\Omega$ i.e. for $\omega\in\Omega$, $\ell (\omega) := \int_0^1 |\dot{\omega}_t|(t) dt$, and let $e_t$ be the time evaluation functional at time $t$.
If $P^*$ is a measure over paths solution to
\begin{equation*}
    \inf_{P \in \mathcal{P}(\Omega)} \int \ell \, d P \text{ s.t. } (e_t)_\# P = \mu_t \text{ for } t \in \{0,1\},
\end{equation*}
then the $W_2$ geodesic between $\mu_0$ and $\mu_1$ is given by $\mu_t = (e_t)_\# P^*$. From this starting point and in the case of $\mathcal{X} = \mathbb{R}^d$, \citet{chen2018measure} define $P$-splines as the minimizers of
\begin{equation}\label{eq:w2_p_spline}
    \inf_{P \in \mathcal{P}(\Omega)} \int c \, d P \text{ s.t. } (e_{t_i})_\# P = \mu_i,
\end{equation}
where this time $\Omega$ is the set of paths with absolutely continuous derivative and $c$ is the curvature cost $c(\omega) = \frac{1}{2} \int_0^1 |\ddot{\omega}_t|^2 \, dt$. Letting $W_2$ $E$-splines be defined by
\begin{equation}\label{eq:w2_e_spline}
    \inf_{(\mu_t,v_t)} \int_0^1 \left \lVert \frac{\mathbf{D}}{dt} v_t \right \rVert_{\mu_t}^2 \, dt \text{ s.t. } \mu_{t_i} = \mu_i,
\end{equation}
\citet{chen2018measure} relates them:
\begin{theorem}\cite[Section 5.3]{chen2018measure}\label{thm:w2_paths_curvature}
    Let $(\mu_t)$ be some curve in $W_2$ with derivative $(v_t)$. Then there is a measure $P \in \mathcal{P}(\Omega)$ such that $(e_t)_\# P = \mu_t$, and the objective in \pref{eq:w2_p_spline} is equal to that of \pref{eq:w2_e_spline}. Furthermore, $P$ can be defined from the flow maps from the continuity equation.
\end{theorem}
This shows that problem \pref{eq:w2_p_spline} is a relaxation of \pref{eq:w2_e_spline}. In \citet[Proposition 2]{chewi2020fast}, it is shown that this is not tight in the sense that there exists (Gaussian) measures $\mu_1,\dots,\mu_n$ that are a solution to problem \eqref{eq:w2_e_spline}, yet are suboptimal for \eqref{eq:w2_p_spline}.

This formulation cannot be directly extended to $\WFR$ since, in \pref{eq:w2_p_spline}, the marginals of $P$ must have the same mass at all times. The solution is to consider curves in a different base space, introduced in \citet{liero2016optimal,chizat2017unbalanced}, which we briefly describe.

\subsection{The Cone Space}

The {\it cone space} $\mathfrak{C}$ is the manifold $\mathbb{R}^d \times \mathbb{R}_+$, with $\mathbb{R}^d \times \{0\}$ identified as a single point.\footnote{Fraktur font characters will be reserved for objects relating to the cone space, consistent with \citet{liero2018optimal}.} Points $(x,r) \in \mathfrak{C}$ are thought of as tuples of position and mass. The distance is given by
\begin{equation*}
    d_\mathfrak{C}\left( (x_0,r_0), (x_1,r_1) \right)^2 = r_0^2 + r_1^2 - 2 r_0 r_1 \cos(|x_0-x_1| \wedge \pi).
\end{equation*}
Call the vertex of the cone (the point $\mathbb{R}^d \times \{0\}$) as $\mathfrak{0}$. Notice that if $|x_0 - x_1| \geq \pi$ then $d_\mathfrak{C}((x_0,r_0),(x_1,r_1)) = r_0 + r_1$. Since for any $(x,r)$ we have $d_\mathfrak{C}((x,r),\mathfrak{0}) = r$, this reflects that fact that if $x_0$ and $x_1$ are far away then the shortest path is to the vertex and back out.

Writing $(v, p) = \frac{d}{dt}(x,r)$, the Riemannian metric on the cone space is given by
\begin{equation}\label{eq:cone_metric}
    \langle (v_1,p_1), (v_2,p_2) \rangle_{(x,r)} = \langle v_1,v_2 \rangle r^2 + p_1 p_2.
\end{equation}

Given a measure $\lambda \in \mathcal{M}(\mathfrak{C})$, we can project it to a measure $\mathfrak{P}\lambda \in \mathcal{M}(\mathbb{R}^d)$ via
\begin{equation*}
    \int f(x) \, d \mathfrak{P}\lambda(x) = \int r^2 f(x) \, d \lambda(x,r).
\end{equation*}
The measure $\lambda$ is then called a lift of $\mathfrak{P} \lambda$. The presence of the term $r^2$, as opposed to $r$, is to simplify the parametrization of $\mathfrak{C}$. Observe that there are many possible liftings of a measure $\mu \in \mathcal{M}(\mathbb{R}^d)$ to $\lambda \in \mathcal{M}(\mathfrak{C})$, the most obvious of which is $d \lambda(x,r) = \delta_1(r) \cdot  d\mu(x)$, where $\delta_p(\cdot)$ stands for the Dirac point measure at $p$. As $\mathfrak{C}$ is (save for the point $\mathfrak{0}$) a Riemannian manifold we can define its $W_2$ metric as usual
\begin{equation*}
    W_{\mathfrak{C},2}^2(\lambda, \eta) = \inf_{\gamma \in \Pi(\lambda, \eta)} \int d_\mathfrak{C}^2 \, d \gamma,
\end{equation*}
where $\Pi(\lambda,\eta)$ is the set of all couplings of measures $\lambda$ and $\eta$ on the cone.
The following useful characterization is proved in \citet{liero2016optimal}.

\begin{theorem}\cite[Theorem 7]{liero2016optimal}\label{thm:wfr_cone_w2}
    For any measures $\mu_0,\mu_1 \in \mathcal{M}_+(\mathbb{R}^d)$, we have
    \begin{equation*}
        \WFR(\mu_0,\mu_1) = \inf_{\lambda_0,\lambda_1}  W_{\mathfrak{C},2}(\lambda_0,\lambda_1),
    \end{equation*}
    where the $\lambda_i \in \mathcal{P}(\mathfrak{C})$ project to $\mu_i$, $\mathfrak{P} \lambda_i = \mu_i$ for $i=0,1$. Furthermore, there are optimal lifts $\lambda_0, \lambda_1$, and an optimal coupling $\gamma$ between them.
\end{theorem}
This allows to characterize WFR geodesics in the same way as in Euclidean space.
\begin{prop}\cite[Section~3.4]{liero2016optimal}\label{thm:wfr_geodeics_cone}
    Let $\mu_0,\mu_1 \in \mathcal{M}_+(\mathbb{R}^d)$, and let $\gamma \in \mathcal{M}(\mathfrak{C}^2)$ be the optimal coupling of the optimal lifts in Theorem \ref{thm:wfr_cone_w2}. For each pair of points $z_0,z_1 \in \mathfrak{C}$, let $g^{z_0,z_1}_t$ be the geodesic between them. Then the curve of measures $(\mu_t)_t$ defined by
    \begin{equation*}
        \mu_t = \mathfrak{P} \left[(g^{z_1,z_2}_t)_\# \gamma \right]
    \end{equation*}
    is a geodesic in $\WFR$ between $\mu_0$ and $\mu_1$.
\end{prop}

We may view the above theorem as saying that $\WFR$ geodesics are identified with certain measures supported on geodesics in $\mathfrak{C}$. In fact, \citet{liero2016optimal} prove that all curves have such a representation, and it is faithful to first order.
\begin{theorem}\cite[Theorem~15]{liero2016optimal}\label{thm:wfr_paths}
    Let $(\mu_t)$ be an absolutely continuous curve in $\WFR$. Then there is a measure $P \in \mathcal{P}(\Omega_\mathfrak{C})$ such that $\mu_t = \mathfrak{P} \left[ (e_t)_\# P \right]$ and for a.e. $t \in [0,1]$
    \begin{equation*}
        |\dot{\mu}_t|^2 = \int |\dot{z}_t|^2 \, d P(z)
    \end{equation*}
    where $|\dot{\mu}_t|$ is the metric derivative, which is equal to the intrinsic Riemannian quantity $\lVert \dot{\mu}_t \rVert_{\mu_t}$.
\end{theorem}
As in the $W_2$ case, this analogy can be extended profitably to second order ---that is, to $|\ddot{\mu}_t|$ in the WFR sense.

\subsection{\texorpdfstring{$P$-}{P}Splines on \texorpdfstring{$\WFR$}{WFR}}

In light of the dual characterization of curvature in Wasserstein space of Theorem \ref{thm:w2_paths_curvature} and the careful reformulation of WFR as $W_2$ on the cone of Proposition \ref{thm:wfr_geodeics_cone} and theorem \ref{thm:wfr_paths}, we define $P$-splines in $\WFR$ by
\begin{equation}\label{eq:wfr_p_spline}
    \inf_{P \in \mathcal{P}(\Omega_\mathfrak{C})} \iint_0^1 |\ddot{z}_t|^2 \, dt \, d P(z) \text{ s.t. } \mathfrak{P} \left[ (e_{t_i})_\# P \right] = \mu_i.
\end{equation}
We prove that this is indeed a relaxation of the $E$-spline problem \pref{eq:wfr_espline}.
\begin{restatable}{theorem}{secondthm}\label{thm:wfr_p_e}
    Let $\mu_t$ be a sufficiently smooth curve in $\WFR$. Then there is a measure $P \in \mathcal{P}(\Omega_\mathfrak{C})$ such that $\mathfrak{P} \left[ (e_t)_\# P \right] = \mu_t$ for all $t$, and the $E$-cost of $\mu$ is equal to the $P$-cost of $P$. The measure $P$ is induced by the flow maps associated to the curve $\mu_t$.
\end{restatable}
\begin{proof}
See appendix \ref{thm:equivalence}.
\end{proof}

To understand the proof (and complete the description of the proposition) we must define the flow maps. In $W_2$, curves of measures are interpretable as particle flows via the {\it flow maps} defined by
\begin{equation}\label{eq:w2_flow_maps}
    \dot{X}_t = v_t(X_t), \, X_0 = \Id.
\end{equation}
It then holds that
\begin{equation*}
    \mu_t = (X_t)_\# \mu_0.
\end{equation*}
There is a similar characterization for sufficiently smooth paths in $\WFR$ due to \citet{maniglia2007probabilistic}.
\begin{prop}\cite[Proposition 3.6]{maniglia2007probabilistic} \label{prop:maniglia_wfr_flow_maps}
    Let $v \in L^1(W^{1,\infty}(\mathbb{R}^d,\mathbb{R}^d),[0,1])$ be a vector field and $\alpha \in \mathcal{C}(\mathbb{R}^d \times [0,1])$ a bounded locally Lipschitz scalar function. For $\mu_0 \in \mathcal{M}_+(\mathbb{R}^d)$, there is a unique weak solution to the nonconservative continuity equation \pref{eq:wfr_ce} with initial measure $\mu_0$. Furthermore, this satisfies
    \begin{equation}
        \mu_t = (X_t)_\#(R_t^2 \cdot \mu_0),
    \end{equation}
    for the flow map $(X_t)$ and scalar field $(R_t)$ which solve the ODE system
    \begin{equation*}
    \begin{cases}
        \dot{X}_t = v_t(X_t), & X_0 = \Id.\\
        \dot{R}_t = 2\alpha_t(X_t)\, R_t, & R_0 = 1.
    \end{cases}
    \end{equation*}
\end{prop}

Analogous to what happens in the Wasserstein case, theorem \ref{thm:wfr_p_e} shows that the spline problem in Lagrangian terms can be relaxed in geometric terms via the flow map characterization of proposition \ref{prop:maniglia_wfr_flow_maps}. Thus, under absolute continuity of the measures and the curve, the Lagrangian and geometric formulations agree up to second order. The first order has three equivalent expressions, namely the metric derivative $|\mu_t'|^2$, the WFR derivative $\|\dot{\mu}_t\|_{\mu_t}^2$ and the Lagrangian form $\int |\dot{z}_t|^2\, dP^*(z)$. A similar equivalence holds for the second order derivatives, i.e. the WFR covariant derivative $\|\nabla_{\dot{\mu}_t}\dot{\mu}_t\|_{\mu_t}^2$ and the Lagrangian form $\int |\ddot{z}_t|^2\,dP^*(z)$. We conjecture that this equivalence holds for higher order derivatives.

\section{Transport Splines}

In this section, we define a tractable and smooth interpolant of measures. Piecewise linear interpolation in $W_2$ has the virtue of being simple to construct, which allows for additional regularization, as it is often found in applications \citep{schiebinger2019optimal,lavenant2021towards}. It yields however a curve that is not smooth. On the other hand, due to the inherent difficulty in solving the $P$-spline problem over $W_2$, \cite{chewi2020fast} define a different, and substantially more tractable, interpolant that they call {\it transport splines}. Very briefly, they proceed as follows:
\begin{enumerate}
    \item Start with measures $\mu_i$ at times $t_i$.
    \item Compute the $W_2$-optimal couplings $\gamma_{i\to i+1}$ from $\mu_i$ to $\mu_{i+1}$, which are induced by maps $T_{i\to i+1}$ provided the $\mu_i$ are absolutely continuous. Let $T_i = T_{0\to 1} \circ \cdots \circ T_{i-1\to i}$
    \item Define for each $x$ the flow map $X_t = S_t[x,T_1x,\ldots,T_Nx]$, where $S_t$ is the natural cubic spline interpolant in Euclidean space (with the base times $t_i$ implicit).
    \item Output $\mu_t = (X_t)_\# \mu_0$.
\end{enumerate}

Not only is computing transport splines very fast, but if the measures $\mu_i$ are sampled from an underlying smooth curve in $W_2$, then the transport spline converges in supremum norm to the true curve at rate $O(\delta^2)$, where $\delta$ is the maximum distance between successive times $t_i$. Importantly for applications, the curve of measure is smooth, which is not true for a piecewise-geodesic interpolant. Other virtues, and their connection with $P$- and $E$-splines, are explored in \citet{chewi2020fast}. In this section we aim to define an analogous interpolant in $\WFR$.

The characterization in Theorem \ref{thm:wfr_cone_w2} is not directly useful here, even if we knew the optimal lifts $\lambda_0$, $\lambda_1$, it is not clear that the optimal coupling is induced by a map. Even if it were, that would not be enough for our purposes --- we would want to associate a unique mass $r$ to each initial point $x$, and map the position-mass pair $(x,r(x))$ to another position-mass pair $(x',r'(x'))$, with $r'$ the unique mass at $x'$. Instead we turn to a third characterization of $\WFR$ given by \citet{liero2018optimal}.

\begin{theorem}(\citealp[Theorem~8]{liero2016optimal};\citealp[Theorem~6.6]{liero2018optimal})\label{thm:wfr_kl}
    Let $\mu_0$ and $\mu_1$ be meausures, and define
    \begin{equation*}
        c(x,y) = -2 \log \cos \left( |x-y| \wedge \tfrac{\pi}{2} \right).
    \end{equation*}
    Then
    \begin{equation}\label{eq:wfr_kl_obj}
        \WFR(\mu_0,\mu_1)^2 = \inf_{\eta \in \mathcal{M}_+} \KL{\eta_0}{\mu_0} + \KL{\eta_1}{\mu_1} + \int c(x,y) \, d \eta(x,y).
    \end{equation}
    Furthermore, the infimum is achieved, and if $\mu_0$ and $\mu_1$ are absolutely continuous with respect to the Lebesgue measure, the optimal $\eta^*$ is unique and is induced by a map. 
\end{theorem}

The cost $c$ differs from the term in the cone metric $d_\mathfrak{C}$ by taking the minimum against $\tfrac{\pi}{2}$, not $\pi$. This is the difference between transport of points in $\mathfrak{C}$ and transport of measures on $\mathfrak{C}$. To transport $z_0 = (x_0,r_0)$ to $z_1 = (x_1,r_1)$ in $\mathfrak{C}$ we may reduce the mass at $x_0$ to $0$, then increase it again at $x_1$ to $r_1$, whereas to transport $\delta_{z_0}$ to $\delta_{z_1}$ these can be done simultaneously (having a superpostion of two deltas), so that the WFR geodesic is at each time a combination of two deltas. The superposition results in a lower overall WFR cost when $|x-y|$ is less than $\frac\pi2$.

The optimal coupling $\eta$ in Theorem \ref{thm:wfr_kl} and the optimal coupling $\gamma$ in Theorem \ref{thm:wfr_cone_w2} are intimately related, as another theorem of \citet{liero2018optimal} shows.
\begin{theorem}\citep[Theorem~6.2]{liero2018optimal}\label{thm:wfr_eta_to_gamma}
    Suppose $\eta$ minimizes the objective in Theorem \ref{thm:wfr_kl}, and let $\eta_i = \pi_i \eta$ be its marginals. Write
    \begin{equation*}
        \mu_i = \sigma_i \eta_i + \mu_i^\bot
    \end{equation*}
    where $\sigma_i = d \eta_i / d \mu_i$ and $\mu_i$ is mutually singular with $\eta_i^\bot$. Define the plan $\gamma_\eta$ by
    \begin{align*}
        d \gamma_\eta(z_0,z_1) &= \delta_{\sqrt{\sigma_0(x_0)}}(r_0) \cdot \delta_{\sqrt{\sigma_1(x_1)}}(r_1) \cdot d \eta(x_0,x_1)\\
        &\quad+ \delta_1(r_0) \cdot d \mu_0^\bot(x_0) \cdot  \delta_\mathfrak{0}(z_1) + \delta_1(r_1) \cdot d \mu_1^\bot(x_1) \cdot \delta_\mathfrak{0}(z_0)
    \end{align*}
    Then $\gamma_\eta$ is optimal for the objective in Theorem \ref{thm:wfr_cone_w2}.
\end{theorem}

Now, suppose the $\mu_0,\mu_1$ are such that there $\mu_i \ll \eta_i$, so that $\mu_i^\bot = 0$ in the theorem above. This happens, for instance, when the conditions of Proposition \ref{prop:maniglia_wfr_flow_maps} are satisfied along the geodesic between them. The optimal $\eta$ for Theorem \ref{thm:wfr_kl} is supported on a map $T$, and thus an optimal $\gamma_\eta$ for Theorem \ref{thm:wfr_cone_w2} is supported on the assignment
\begin{equation*}
    \left(x_0, r_0(x_0)\right) \to \left(T(x_0), r_1(T(x_0))\right)
\end{equation*}
which associates to each $x_0$ a unique mass $r_0(x_0) = \sqrt{\sigma_0(x_0)}$ and maps it to another unique location and mass $r_1(T(x_0)) = \sqrt{\sigma_1(T(x_0))}$. This is much stronger than $\gamma_\eta$ being induced by a map on $\mathfrak{C}$.

Now, define the operator $S^\mathfrak{C}_t[z_0,\ldots,z_N]$ to be the Riemannian cubic interpolant in $\mathfrak{C}$ of the points $z_0,\ldots, z_N$. We define {\it transport splines over} $\WFR$ by the following procedure
\begin{enumerate}
    \item For each $i$ solve \pref{eq:wfr_kl_obj} between $\mu_i$ and $\mu_{i+1}$ to obtain the optimal coupling $\eta_{i\to i+1}$, and thus the map $T_{i\to i+1}$. Let $T_i = T_{0\to 1} \circ \cdots \circ T_{i-1 \to i}$.
    \item For each $x$, form a path $\widetilde{X}_t(x)$ interpolating the $x_i = T_i(x)$ and a mass path $\widetilde{R}_t(x)$ interpolating the masses $r_i(x_i)$.
    \item Define the interpolating curve on $[t_{i-1},t_i]$ by computing the Cone spline $(X_t,R_t)$ with endpoint velocity constraints $\dot{X}_{t_j}(x) = \frac{d}{dt} \widetilde{X}_{t}(x)|_{t=t_j}$ and $\dot{R}_{t_j}(x) = \frac{d}{dt} \widetilde{R}_t(x)|_{t=t_j}$ for $j=i-1,i$ and $i=1,\dots, n$.
\end{enumerate}
In the first step, we obtain a family of transport maps and corresponding masses at each point $x$ from the coupling $\eta_{i\to i+1}$. Notice that two adjacent couplings $\eta_{i-1\to i}$ and $\eta_{i\to i+1}$ may not have the same $i$-th marginal, so the mass $\sqrt{\sigma_i}$ may be not uniquely defined. We remedy this by rescaling the coupling to have unit mass, so that the coupling is supported on $\sqrt{\mu_i(x_i)}$. Lifted plans on the cone are scale-invariant, thus the optimality of \eqref{eq:wfr_kl_obj} is preserved under this scaling, as explained in \citet[Section~\S3.3]{liero2016optimal}. In the next step, we would ideally interpolate any sequence $(x_i,r_i)_i$ using Riemannian cubics on $\mathfrak{C}$, these are difficult to compute; unlike the Euclidean case, there appears to be no closed formula. We propose instead approximating the velocities at each knot from the (Euclidean) spline curves of $X$ and $R$ independently \citep[Appendix D]{chewi2020fast}. We then solve in step 3 the cone spline problem by specifying the velocities computed in the previous step and using De Casteljau's algorithm.

\subsection{De Casteljau's algorithm on the cone}

In Euclidean space, an interpolating curve with minimal curvature is completely determined by the endpoint velocities. The De Casteljau's algorithm computes points on this curve by iteratively finding points along the paths between the endpoints and two other control points. We describe the procedure for reference. The De Casteljau curve that interpolates $x_0$ and $x_3$ with speeds $v_0$ and $v_3$ at times $t=0$ and $t=1$ respectively is constructed as follows:
\begin{enumerate}
    \item Compute control points $x_1=x_0 + \frac{v_0}{3}$ and $x_2 = x_2 - \frac{v_3}{3}$.
    \item Compute the first intermediate points $w_i(t) = (1-t)x_i + t x_{i+1}$, for $i=0,1,2$.
    \item Compute the second intermediate points $u_j(t) = (1-t)w_j + t w_{j+1}$, for $j=0,1$.
    \item Compute the spline as  $p(t) = (1-t)u_0 + t u_1$.
\end{enumerate}

We extend this definition to the cone by replacing the linear interpolation by a geodesic interpolation \citep{Absil2016,Absil2019}. It is known that the geodesic has closed-form expression given by:

\begin{theorem}\citep[Section~\S8.1]{liero2018optimal}\label{thm:geodesics}
The geodesic interpolator of the points $z_j = (x_j,r_j)$ is given by $z\left(t \right)=\left(x\left(t \right), r\left(t \right)\right)$, where
\[ 
\begin{aligned}
x\left(t \right) &=\left(1-\rho\left(t \right)\right) x_{0}+\rho\left(t \right) x_{1},\\
r\left(t \right)^{2} &=(1-t)^{2} r_{0}^{2}+t^{2} r_{1}^{2}+2 t(1-t) r_{0} r_{1} \cos \left|x_{0}-x_{1}\right|, \\
\rho(t) &= \frac{1}{|x_1-x_0|} \arccos\left( \frac{(1-t)r_0 + t r_1\cos|x_1-x_0|}{r(t)}\right).
\end{aligned}
\]
\end{theorem}

Theorem \ref{thm:geodesics} is only valid for distances $|x_1-x_0| < \pi$, otherwise the geodesic simply goes from $x_0$ through the tip of the cone $\mathfrak{o}$ and back to $x_1$. In terms of measures, this behavior corresponds to pure growth-decay and no transport. From the discussion above about geodesics on the cone and geodesics on WFR, and because we want to model scenarios where both transport and growth are present at every moment, we restrict the computation of cone splines to points that are contained within a ball of radius $\pi/2$. We will denote the geodesic interpolation of points $z_0$ and $z_1$ at time $t$ as $z_0\#_t z_1$.

In order to extend step 1 above, we must first understand how the velocities of a De Casteljau spline at $t=0,1$ relate to the control points $x_1,x_2$.

\begin{restatable}{prop}{thirdthm}\label{prop:cone_decasteljau}

Let $(x_0,r_{x_0})$ and $(x_3,r_{x_3})$ be given points on the cone, and $(x_1,r_{x_1})$ and $(x_2,r_{x_2})$ be control points. Let 
\begin{align*}
 (w_i(t),r_{w_i}(t)) &:= (x_i,r_{x_i}) \#_t  (x_{i+1},r_{x_{i+1}}),\quad i = 0,1,2  \\ 
 (u_j(t),r_{u_j}(t)) &:= (w_j,r_{w_j}) \#_t  (w_{j+1},r_{w_{j+1}}), \quad j = 0,1   \\
 (p(t),r_p(r)) &:= (u_0,r_{u_0}) \#_t  (u_{1},r_{u_{1}}),   
\end{align*}
be the De Casteljau spline on the cone, then 
\begin{align} 
\dot{p}(0) &= 3 \frac{r_{x_1}}{r_{x_0}} \frac{\sin|x_0-x_1|}{|x_0-x_1|} (x_1-x_0) \label{eq:endpoints_p0} \\
\dot{r}_p(0) &= 3 (r_{x_1} \cos |x_0-x_1| - r_{x_0}) \label{eq:endpoints_r0} \\
\dot{p}(1) &= 3 \frac{ r_{x_2}}{r_{x_3}} \frac{\sin|x_3-x_2|}{|x_3-x_2|} (x_{3}-x_2) \label{eq:endpoints_p1},\\
\dot r_{p}(1) &= 3 \left( r_{x_{3}} - r_{x_2}\cos|x_3-x_2| \right) \label{eq:endpoints_r1}.
\end{align}
\end{restatable}
\begin{proof}
See Appendix \ref{thm:cone_de_casteljau}.
\end{proof}

In particular, at an interval $[t_i,t_{i+1}]$, in order to achieve the position and mass velocities given by $v_i=\frac{d}{dt} \widetilde{X}_{t}(x)|_{t=t_i}$ and $s_i = \frac{d}{dt} \widetilde{R}_{t}(x)|_{t=t_i}$, we have to choose the control point $(x_1,r_1)$ as $x_1 = x_{t_i} + c_1 \frac{v_i}{\|v_i\|}$ and $r_1 = c_2 r_{t_i}$, where $c_1,c_2$ are solved from (\ref{eq:endpoints_p0}-\ref{eq:endpoints_r0}):
\begin{align*}
c_1 &= \arctan\left(\frac{\|v_i\|}{\tfrac{s_i}{r_{t_i}} +3 \delta^{-1}}\right),  \\
c_2 &=  \tfrac{\delta}{3} \sqrt{\|v_i\|^2+\left(\tfrac{s_i}{r_{t_i}}+3 \delta^{-1}\right)^2},
\end{align*}
where $\delta = t_{i+1} - t_i$.  Similarly for $t_{i+1}$, we choose $(x_2,r_2)$ with $x_2 = x_{t_{i+1}} -  c_3 \tfrac{v_{i+1}}{\|v_{i+1}\|}$, $r_2 = c_4 r_{t_{i+1}}$, and solve (\ref{eq:endpoints_p1}-\ref{eq:endpoints_r1}):
\begin{align*}
c_3 &= \arctan\left(\frac{\|v_{i+1}\|}{3\delta^{-1} - \tfrac{s_{i+1}}{r_{t_{i+1}}}} \right),  \\
c_4 &=  \tfrac{\delta}{3} \sqrt{\|v_{i+1}\|^2+\left(3\delta^{-1} - \tfrac{s_{i+1}}{r_{t_{i+1}}} \right)^2}.
\end{align*}

The derivation of the above equations requires that $c_1,c_3 \ge 0$; it reflects the fact that $x_1-x_0$ and $x_3-x_2$ point in the right direction. This translates to $s_i \ge -3 \delta^{-1} r_{t_i} $ and $s_{i+1} \le 3 \delta^{-1} r_{t_{i+1}}$, which constraints the possible velocities that are achievable by the De Casteljau's algorithm. To comprehend why this happens, notice that unlike the Euclidean case, the second variable $r$ is constrained to be non-negative. This limits the possible masses of the control points $x_1,x_2$, which in turn limits the endpoint derivatives $\dot{r}_p(0)$ and $\dot{r}_p(1)$ (equations \eqref{eq:endpoints_r0},\eqref{eq:endpoints_r1}). 

There are two limitations to the cone transport spline problem: the maximum $\frac\pi2$ diameter in the position space and the bounds in knot velocities in mass space, as described above. Furthermore, the velocities are estimated from two independent spline problems. A simple solution to these limitations is to 1) scale down all distances to a space with a smaller diameter and 2) scale up the knot times to allow for smaller endpoint derivatives.

\section{Numerical Experiments}

In this section, we illustrate the behavior of the cone transport spline via the De Casteljau's algorithm. We showcase different aspects of the paths qualitatively by considering spline interpolation of measures in one and two dimensions. In section \ref{sec:time} we consider interpolation with varying times. We show in section \ref{sec:space} two alternative ways of solving a spline problem given by discrete measures by considering a problem in two dimensions.\footnote{The MATLAB code used for generating the figures can be found here \url{https://github.com/felipesua/WFR_splines}.} %TODO Finally, we apply the cone transport spline algorithm to the single-cell dataset [Expand] in \parencite{schiebinger2019optimal} in section XX.
% and compare it to linear interpolation with successive transport maps

\subsection{Time effects and Linear interpolation}\label{sec:time}

\begin{figure}
    \centering
    \includegraphics[width=.3\textwidth]{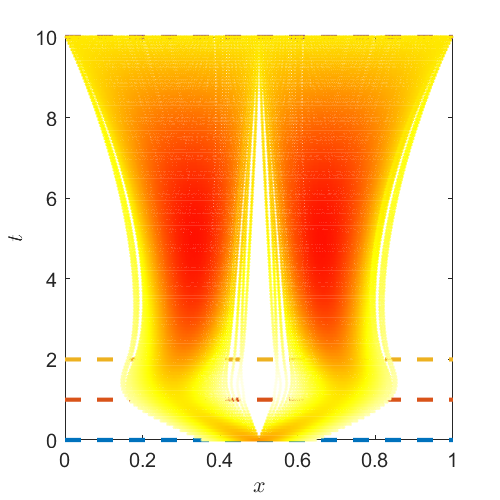}
    \includegraphics[width=.3\textwidth]{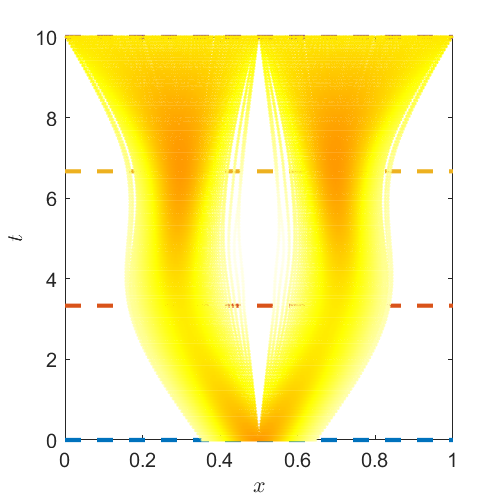}
    \includegraphics[width=.3\textwidth]{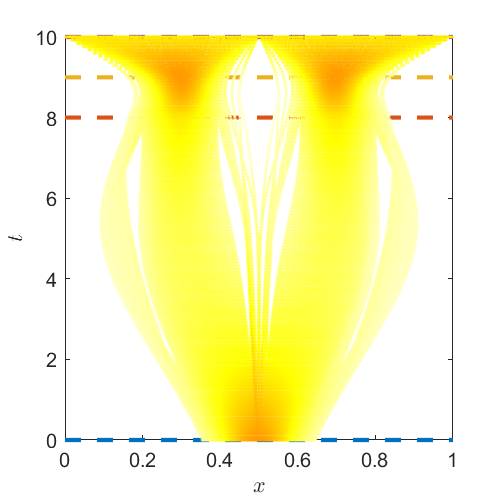}
    
    \includegraphics[width=.3\textwidth]{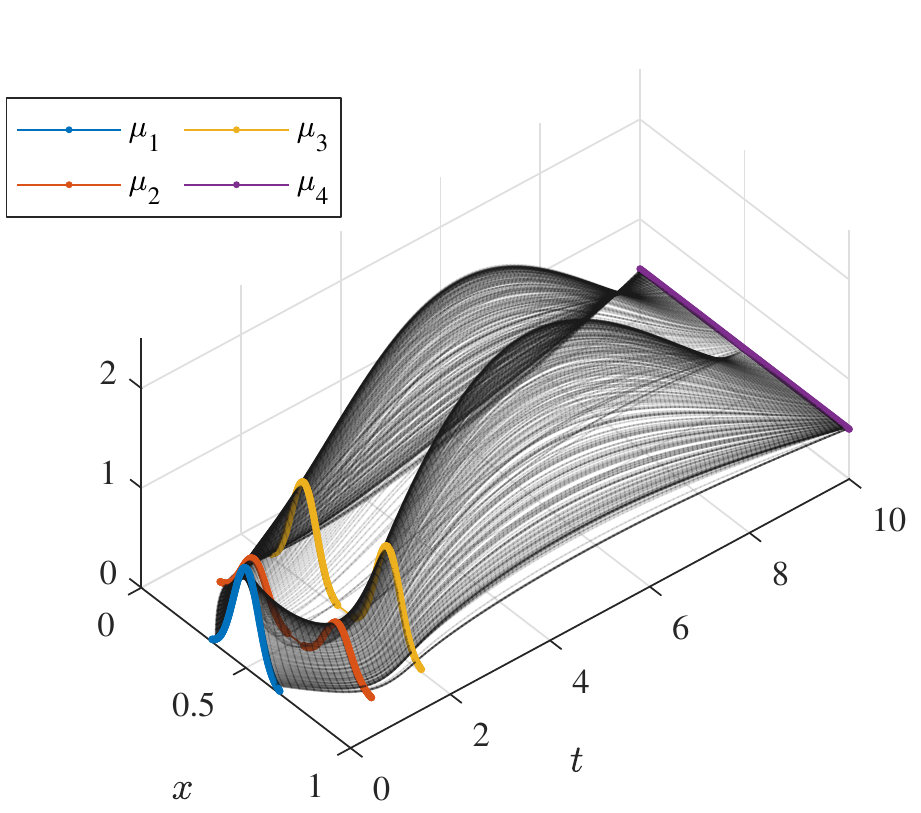}
    \includegraphics[width=.3\textwidth]{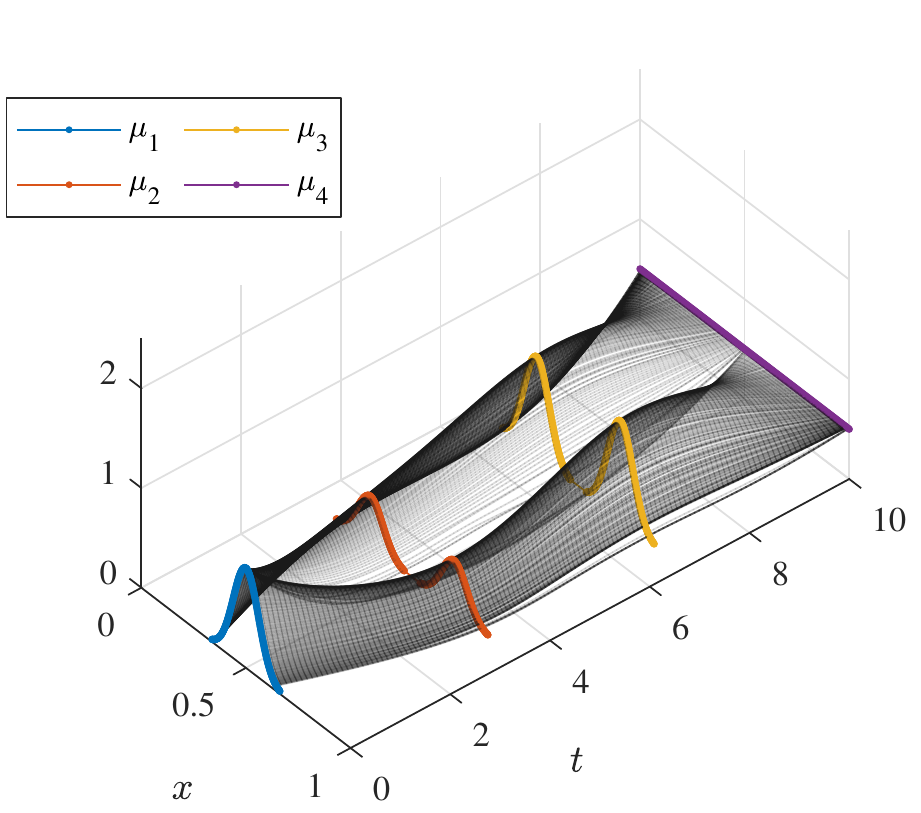}
    \includegraphics[width=.3\textwidth]{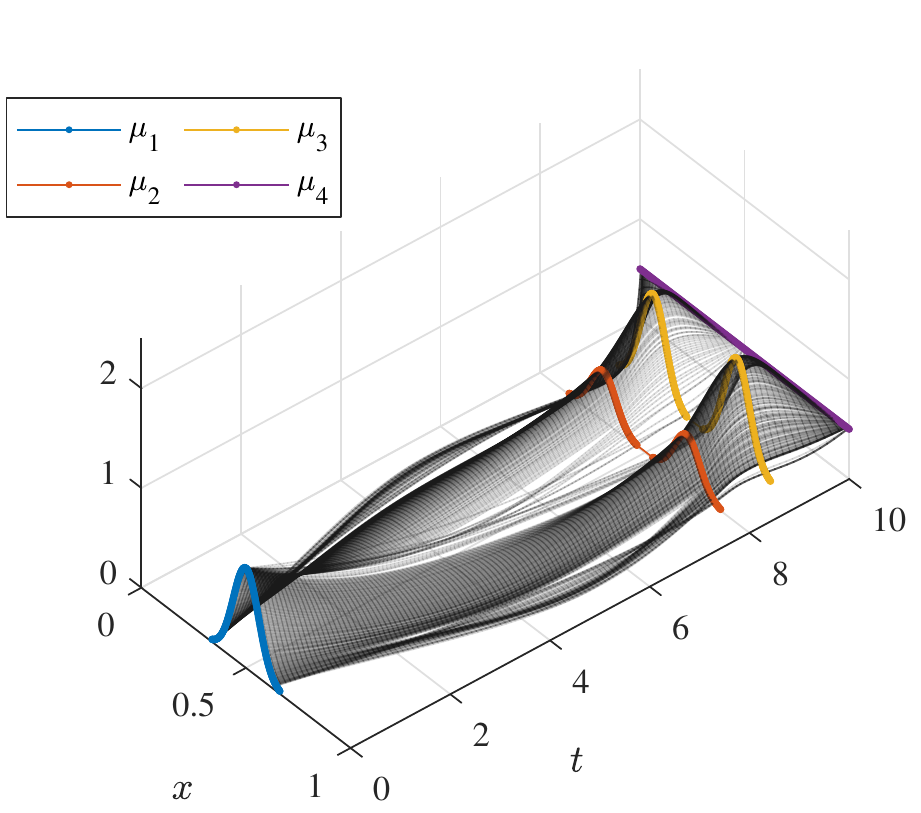}
    \caption{\footnotesize De Casteljau cone interpolation of measures $\mu_1,\dots,\mu_4$. Measure values are shown as a heat map (top) and in a three-dimensional plot (bottom) as a function of time. Knot times $(t_1,\dots,t_4)$ are $(0,1,2,10)$ (left), $(0,\tfrac{10}{3},\tfrac{20}{3}, 10)$ (center), and $(0,8,9,10)$ (right). $\sigma=0.06$.}
    \label{fig:one_dim}
\end{figure}
 
In this subsection, we look at the effects of time distribution of knots for a one dimensional cone spline problem. The map obtained from solving \eqref{thm:wfr_kl} is induced by a map only when the measures are absolutely continuous with respect to the Lebesgue measure. Indeed, the WFR plan from a single particle to two particles cannot be induced by a map. The same effect is present when solving problem \eqref{thm:wfr_kl} computationally. Here we choose to solve problem \eqref{thm:wfr_kl} using entropic regularization \citep{chizat2018scaling,peyre2019computational} and derive an approximation of the map from its solution by taking expectation\footnote{We can assume $\eta$ is a probability measure by the scale invariance property of the projection from the cone. } of the marginals: $T(x) = \mathbb{E}_\eta[y|x]$, inspired by \citet{pooladian2021entropic}.

Let $\gamma_\sigma(x) = \exp(-\tfrac{1}{2\sigma^2}x^2) \mathds{1}_{[-2,2]}(x/\sigma).$ For our first experiment, we interpolate the measures $$d\mu_1(x) = \gamma_\sigma\left(x-\tfrac12\right)\,dx,$$ 
$$d\mu_{2}(x) = \tfrac12 \left( \gamma_\sigma(x-0.3) + \gamma_\sigma(x-0.7)\right)\,dx,$$ 
$$d\mu_{3}(x) = \left( \gamma_\sigma(x-0.3) + \gamma_\sigma(x-0.7)\right)\,dx, $$
$$d\mu_{4}(x) =  \tfrac12 \mathds{1}_{[0,1]}(x)\,dx, $$
which we interpolate at three set of knot times $(t_1,\dots,t_4)$: $(0,1,2,10)$, $(0,\tfrac{10}{3},\tfrac{20}{3}, 10)$, and $(0,8,9,10)$ in figure \ref{fig:one_dim}.

\subsection{Space discretization} \label{sec:space}

\begin{figure}
    \centering
    \includegraphics[width=.45\textwidth]{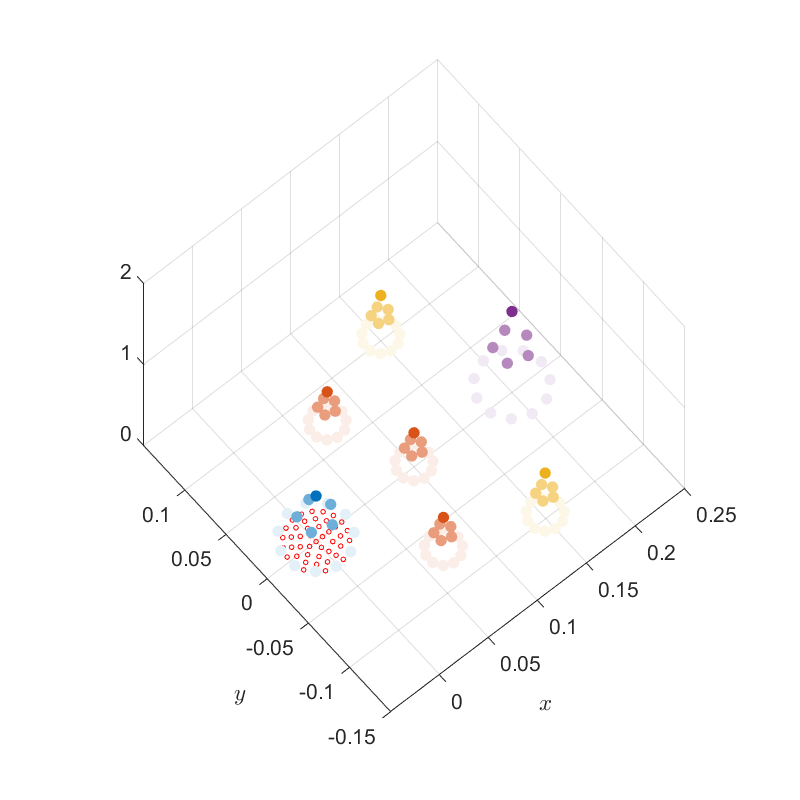}
    \includegraphics[width=.45\textwidth]{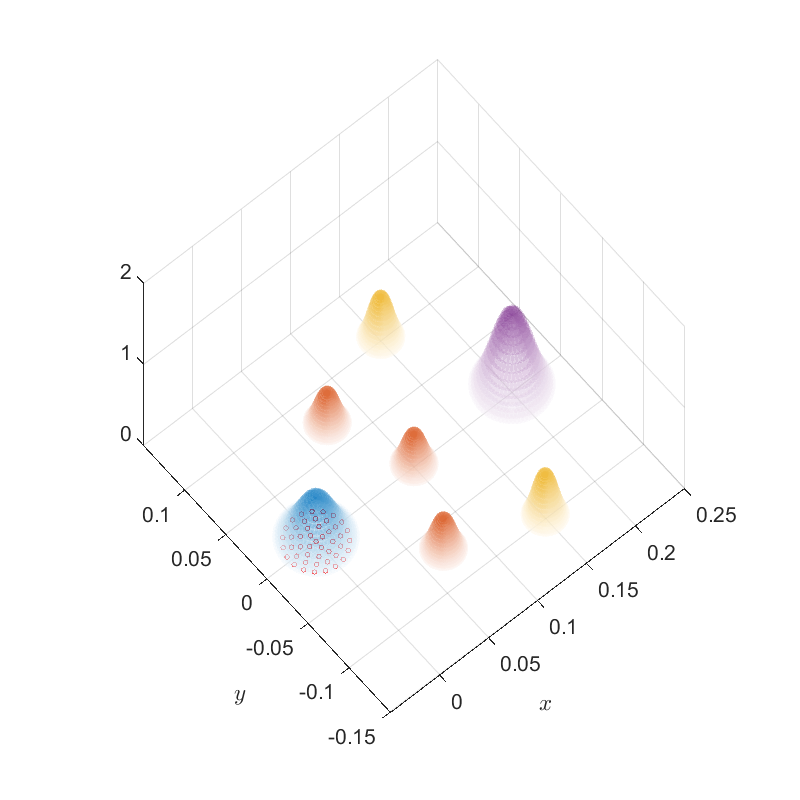}
    
    \includegraphics[width=.45\textwidth]{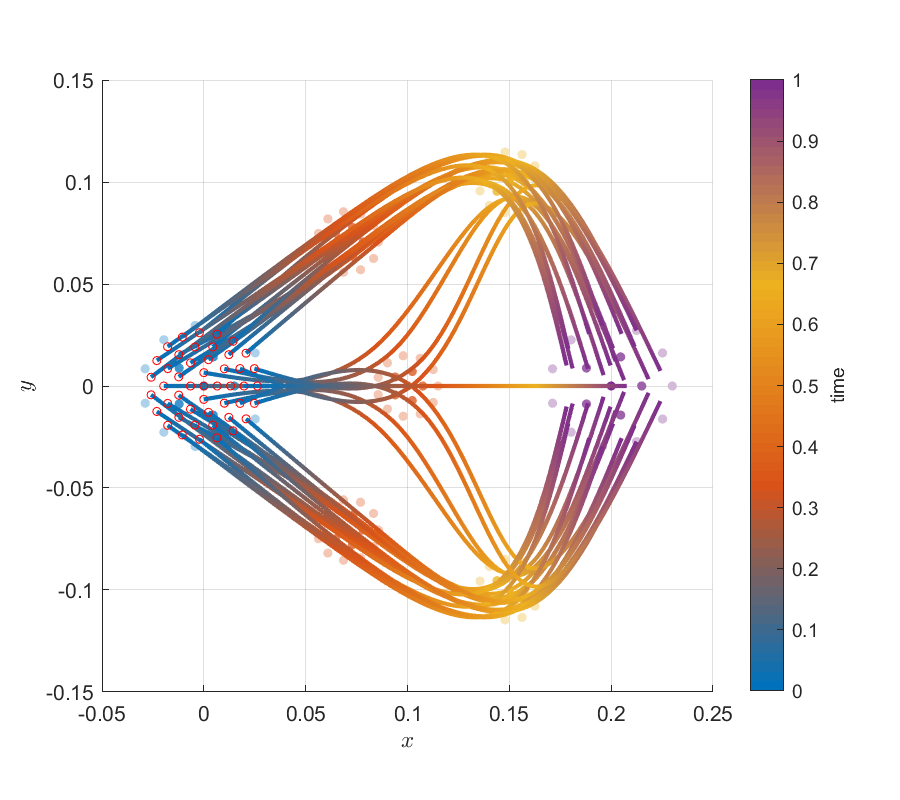}
    \includegraphics[width=.45\textwidth]{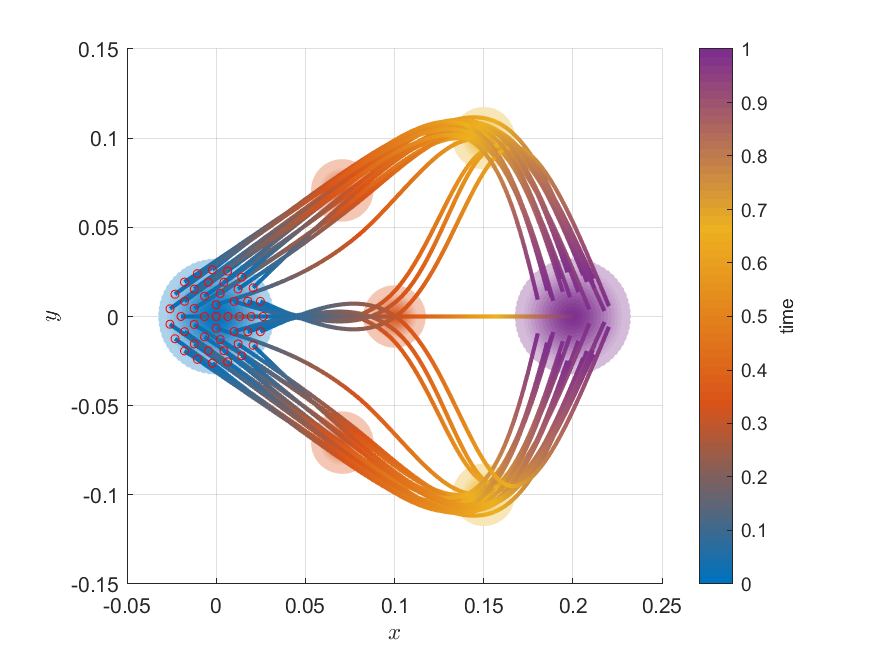}
    
    \includegraphics[width=.45\textwidth]{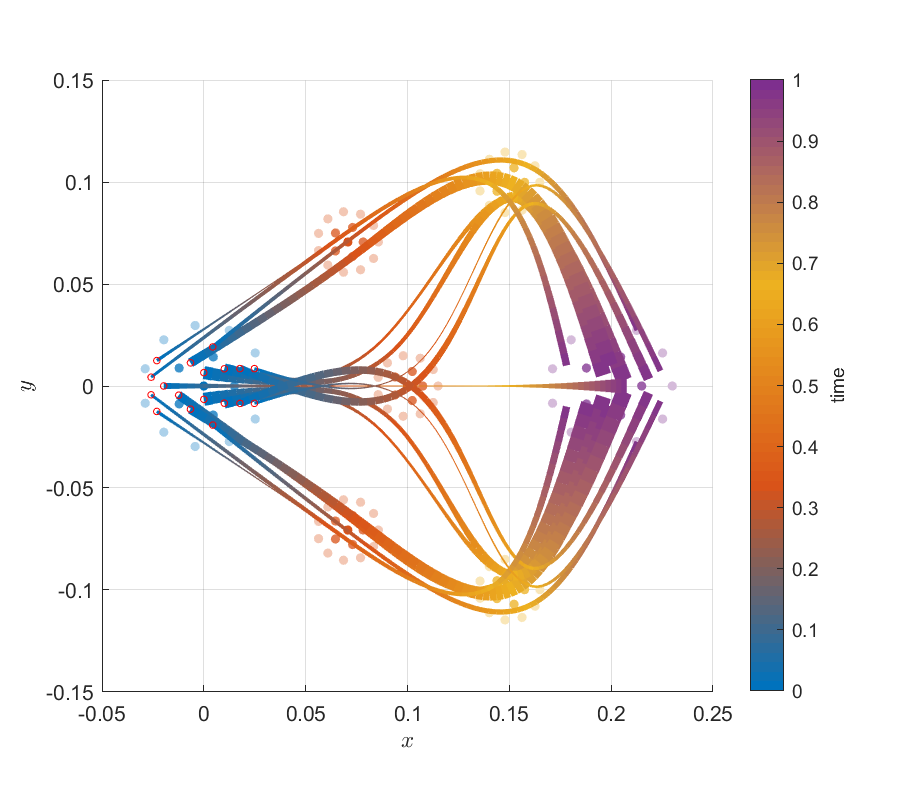}
    \includegraphics[width=.45\textwidth]{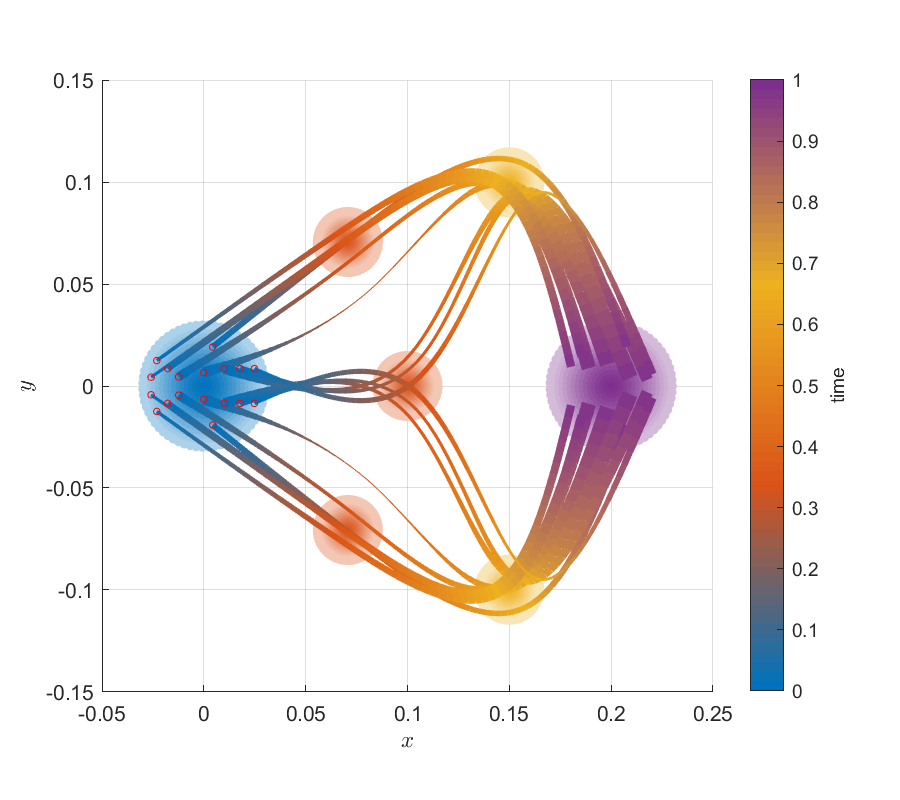}
    \caption{\footnotesize Subsampled (left) and gridded (right) interpolating measures. Measures $\mu_1,\dots,\mu_4$ are shown in blue, orange, yellow, and purple respectively (top). Interpolated paths are computed for the points shown in red (middle) as a function of time varying by color. 
    The width of select curves is changed according to the value of the measure (bottom). $\sigma=0.01$}
    \label{fig:two_dim}
\end{figure}

We give details for the interpolation of measures in two dimensions in figure \ref{fig:two_dim}. For this example suppose we want to interpolate a set of particles. Since we need to represent them as absolutely continuous measures to be able to obtain a map, we convolve with a kernel and discretize the domain like in the previous section. As noted by \citet{lavenant2021towards}, gridding grows exponentially with the dimension which makes it infeasible for high-dimensional applications. We show here the result of fine gridding and compare it to sampling uniformly from the support.

Suppose the resulting interpolating measures are as follows. Let $$\gamma_\sigma(x) = \exp(-\tfrac{1}{2\sigma^2}\|x\|^2) \mathds{1}_{B_2(0)}(x).$$ We interpolate 

$$ d\mu_1(x) = \tfrac34 \gamma_{2\sigma}(x)\, dx, $$
$$ d\mu_2(x) = 0.65 \left( \gamma_{\sigma}(x-(\tfrac{\sqrt{2}}{20},\tfrac{\sqrt{2}}{20})) + \gamma_{\sigma}(x-(0,-\tfrac{\sqrt{2}}{20})) + \gamma_{\sigma}(x-(\tfrac{\sqrt{2}}{20},-\tfrac{\sqrt{2}}{20})) \right) \, dx, $$
$$ d\mu_2(x) = \tfrac34 \left( \gamma_{\sigma}(x-(\tfrac{3}{20},\tfrac{3}{20})) + \gamma_{\sigma}(x-(\tfrac{3}{20},-\tfrac{3}{20})) \right) \, dx, $$
$$ d\mu_4(x) = \gamma_{2\sigma}(x-(\tfrac15,0))\, dx. $$

\section*{Acknowledgements}
We thank Philippe Rigollet and Sinho Chewi for helpful comments and suggestions on the manuscript.
\appendix 

\section{Derivation of the covariant derivative}

\firstthm*
 \begin{proof} \label{thm:covariant_derivative}
  Let $\mathbf{u}_t^i = (u_t^i,\beta_t^i)$ be two tangent fields along a curve $\mu_t$, which has derivative $(v_t,\alpha_t)$. Metric compatibility then reads
\begin{equation*}\label{eq:metric_comp_wfr_1}
\begin{aligned}
    \frac{d}{dt} \langle \mathbf{u}_t^1, \mathbf{u}_t^2 \rangle_{\mu_t} &= \int \langle \partial_t u_t^1, u_t^2 \rangle + \langle u_t^1, \partial_t u_t^2 \rangle + 4 \partial_t \beta_t^1 \, \beta_t^2 + 4 \beta_t^1 \, \partial \beta_t^2 \, d \mu_t\\
    &\quad + \int \langle u_t^1, u_t^2 \rangle + 4 \beta_t^1 \beta_t^2 \, d (\partial_t \mu_t)\\
    &= T_1 + T_2
\end{aligned}
\end{equation*}
By the continuity equation and the dual definition of $\Div v_t \mu_t$, the second integral becomes
\begin{align*}
    T_2 &= \int \langle \nabla u_t^1 \cdot v_t, u_t^2 \rangle + \langle u_t^2, \nabla u_t^2 \cdot v_t \rangle + 4 \alpha_t \langle u_t^1, u_t^2 \rangle\\
    &\quad + \int 4 \beta_t^2 \langle \nabla \beta_t^1, v_t \rangle + 4 \beta_t^1 \langle \nabla \beta_t^2, v_t \rangle + 16 \alpha_t \beta_t^1 \beta_t^2 \, d \mu_t
\end{align*}
Now, we need this to be equal to $\langle \mathbf{u}_t^1, \frac{\mathbf{D}}{dt} \mathbf{u}_t^2 \rangle + \langle \frac{\mathbf{D}}{dt}  \mathbf{u}_t^1,  \mathbf{u}_t^2 \rangle$. Some terms are uniquely attributable to one or the other, such as $\int \langle \partial_t u_t^1, u_t^2 \rangle$, but some are not, such as $4 \beta_t^2 \langle \nabla \beta_t^1, v_t \rangle$. This can arise from either of the two terms: as $\nabla \beta_t^i = u_t^i$,
\begin{equation*}
    \left \langle (0, \langle \nabla\beta_t^1, v_t \rangle), \mathbf{u}_t^2 \right \rangle_{\mu_t} = \left \langle \mathbf{u}_t^1, (4 \beta_t^2 v, 0) \right \rangle_{\mu_t}
\end{equation*}
Indeed, by splitting these terms and gathering the others, the possible covariant derivatives that satisfy metric compatibility are of the form
\begin{equation*}
    \frac{\mathbf{D}}{dt} \mathbf{u}_t = \mathcal{P}_{\mu_t} \begin{pmatrix}
    \partial_t u_t + \nabla u_t \cdot v_t + 2 \alpha_t u_t + 4 p \beta_t v_t\\
    \partial_t \beta_t + (1-p) \langle \nabla \beta_t, v_t \rangle + 2 \alpha_t \beta_t
    \end{pmatrix}
\end{equation*}
for real $p$. Let us check the torsion-free identity with $p=1/2$. In this case, if $\mathbf{u}_t = (\nabla \phi, \phi)$ is constant in time, then
\begin{equation*}
    \frac{\mathbf{D}}{dt} \begin{pmatrix}\nabla \phi\\ \phi\end{pmatrix} = \mathcal{P}_{\mu_t}
    \begin{pmatrix}
        \nabla^2 \phi \cdot v_t + 2 \alpha_t \nabla \phi + 2 \phi v_t\\
        \frac{1}{2}\langle \nabla \phi, v_t \rangle + 2 \alpha_t \phi
    \end{pmatrix}
\end{equation*}
Now, with the setup as in the $W_2$ case, defining $F \colon \phi \mapsto \int \phi \, d \mu$ and writing $\bm{\phi} = (\nabla \phi, \phi)$, we have from the continuity equation
\begin{equation*}
    \partial_t F[\mu_t^i] = \int \langle \bm{\phi}, \mathbf{v}_t^i \rangle \, d \mu_t^i
\end{equation*}
where $\mathbf{v}_t^i$ is the derivative of $\mu_t^i$. As above, we have
\begin{align*}
    \mathbf{u}_0^1(\mathbf{u}^2(F))[\mu] &= \frac{d}{dt} \langle \bm{\phi}, \mathbf{u}_t^2 \rangle_{\mu_t^2} \big|_{t=0}\\
    &= \left \langle \frac{\mathbf{D}}{dt} \bm{\phi}, \mathbf{u}_t^2 \right \rangle_{\mu_t^2} + \left \langle \bm{\phi}, \nabla_{\mathbf{u}_0^1} \mathbf{u}_t^2 \right \rangle_{\mu_t^2} \bigg|_{t=0}
\end{align*}
Recalling that at $t=0$ we have $\mathbf{u}_0^2 = \mathbf{v}_0^1 = (v_0^1, \alpha_0^1)$, the first term becomes (we may ignore the projection, since $\mathbf{u}_t^2$ is already tangent)
\begin{equation*}
    Q_1 = \left \langle \begin{pmatrix}
        \nabla^2 \phi \cdot v_0^2 + 2 \alpha_0^2 \nabla \phi + 2 \phi v_0^2\\
        \frac{1}{2} \langle \nabla \phi, v_0^2 \rangle + 2 \alpha_0^2 \phi
    \end{pmatrix}, \begin{pmatrix}
        v_0^1\\
        \alpha_0^1
    \end{pmatrix} \right \rangle_\mu
\end{equation*}
while the corresponding term from $\mathbf{u}_0^2(\mathbf{u}^2(F))[\mu]$ is
\begin{equation*}
    Q_2 = \left \langle \begin{pmatrix}
        \nabla^2 \phi \cdot v_0^1 + 2 \alpha_0^1 \nabla \phi + 2 \phi v_0^1\\
        \frac{1}{2} \langle \nabla \phi, v_0^1 \rangle + 2 \alpha_0^1 \phi
    \end{pmatrix}, \begin{pmatrix}
        v_0^2\\
        \alpha_0^2
    \end{pmatrix} \right \rangle_\mu
\end{equation*}
and we must check that these agree. The ``top-left'', ``top-right'', and ``bottom-right'' terms are identical for both. The top-middle of the first is equal to the bottom-left for the second, and vice-versa. Thus we have shown that the covariant derivative is given by
\begin{equation*}
    \frac{\mathbf{D}}{dt} \mathbf{u}_t = \mathcal{P}_{\mu_t} \begin{pmatrix}
    \partial_t u_t + \nabla u_t \cdot v_t + 2 \alpha_t u_t + 2 \beta_t v_t\\
    \partial_t \beta_t + \frac{1}{2}\langle \nabla \beta_t, v_t \rangle + 2 \alpha_t \beta_t
    \end{pmatrix}
\end{equation*}
and in specific,
\begin{equation}
    \frac{\mathbf{D}^2}{dt^2} \mu_t = \begin{pmatrix}
        \partial_t v_t + \nabla v_t \cdot v_t + 4 \alpha_t v_t\\
        \partial_t \alpha_t + \frac{1}{2}|\nabla \alpha_t|^2 + 2 \alpha_t^2
    \end{pmatrix}
\end{equation}
This quantity is tangent, so no projection is necessary. \qedhere
\end{proof}

\section{Proof of equivalence in curvature}

\secondthm*
\begin{proof}\label{thm:equivalence}
    Making the definition explicit, we wish of our measure $P$ that the cost of \pref{eq:wfr_p_spline} is equal to
    \begin{equation}\label{eq:wfr_e_cost_explicit}
        \int_0^1 \int \left| \partial_t v_t + \nabla v_t \cdot v_t + 4 \alpha_t v_t \right|^2 + 4 \left( \partial_t \alpha_t + \tfrac{1}{2}|\nabla \alpha_t|^2 + 2 \alpha_t^2 \right)^2 \, d \mu_t \, dt
    \end{equation}
    Let $\lambda_0 \in \mathcal{P}(\mathfrak{C})$ be a lift of $\mu_0$. By Proposition \ref{prop:maniglia_wfr_flow_maps}, the measure $\lambda_t = (X_t, R_t)_\# \lambda_0$ is a lift of $\mu_t$. 
    
    In order to compute the covariant derivative of a curve on the cone, we compute the Christoffel symbols. These are given by the formulas
    \begin{equation*}
        \nabla_{\partial x_i}\partial x_j = \Gamma_{ij}^k \partial_k, \quad \Gamma_{ij}^k = \frac{1}{2}g^{kl}\left(\frac{\partial}{\partial_j} g_{il} + \frac{\partial}{\partial_i} g_{jl} - \frac{\partial}{\partial_l} g_{ij}\right).
    \end{equation*}
    Let $\partial x_1,\dots,\partial x_n,\partial _r$ be the coordinate basis of the tangent space on the cone, hence
    \begin{equation*}
        (g_{ij})_{ij} = \begin{pmatrix} r^2 I_n & 0\\ 0 & 1 \end{pmatrix}, \quad \frac{\partial}{\partial_k}g_{ij} = \begin{cases} 2r & i=j\neq r,~k=r, \\ 0, & \text{else}.\end{cases}
    \end{equation*}
    Since $(g^{ij})$ and $\partial_r g_{ij}$ are only defined along the diagonal and $\partial_l g_{ij}=0$ for $l\neq r$, the only terms that do not vanish are 
    \begin{equation*}
        \Gamma_{rj}^k = r^{-1}\delta_{jk}, \quad \Gamma_{ij}^r = -r\delta_{ij}
    \end{equation*}
    thus the Levi-Civita connection on the cone is given by 
    \begin{equation*}
    \nabla_{\partial_r} \partial_r = 0, \quad\quad \nabla_{\partial_r} X = r^{-1} X, \quad\quad \nabla_{X_1}X_2 = -r\langle X_1,X_2\rangle \partial_r.
    \end{equation*}
    
    From this, the covariant derivative of a curve $z_t = (x_t, r_t)$ on the cone is given by
    \begin{equation*}\label{eq:cone_d2dt2}
        \ddot{z}_t = \nabla_{\dot{z}_t}\dot{z}_t = \left(\ddot{x}_t + 2 \frac{\dot{r}_t}{r_t} \dot{x}_t, ~\ddot{r}_t - r_t |\dot{x}_t|^2 \right)
    \end{equation*}

    thus from the Riemannian metric \pref{eq:cone_metric}
    \begin{equation*}
        |\ddot{z}_t|^2 = \left|r_t \ddot{x}_t + 2 \dot{r}_t \dot{x}_t  \right|^2 + \left|\ddot{r}_t - r_t |\dot{x}_t|^2 \right|^2
    \end{equation*}
    Let $\lambda_0$ be any lift of $\mu_0$ and define $P$ to place mass $\lambda_0(z_0)$ on the path $z_t = (X_t(x_0),R_t(r_0))$, so that $P$ is supported on the flow map curves in $\Omega_\mathfrak{C}$. By applying the total derivative to the defining equations of the flow maps, these curves satisfy
    \begin{align*}
        \dot{x}_t &= v_t(x_t)\\
        \ddot{x}_t &= \partial_t v_t + \nabla v_t \cdot v_t\\
        \dot{r}_t &= 2 r_t \alpha_t(x_t)\\
        \ddot{r}_t &= 2 r_t \left( \partial_t \alpha_t + \nabla \alpha_t \cdot v_t + 2 \alpha_t^2 \right)
    \end{align*}
    Now, we have
    \begin{equation*}
        \int_{\Omega_\mathfrak{C}} \int_0^1 |\ddot{z}_t|^2 \, dt \, dP(z) = \int_\mathfrak{C} \int_0^1 |\ddot{z}_t|^2 \, dt \, d \lambda_0(z_0)
    \end{equation*}
    Expanding out, this is
    \begin{equation*}
        \int_\mathfrak{C} \int_0^1 \left|r_t \ddot{x}_t + 2 \dot{r}_t \dot{x}_t  \right|^2 + \left|\ddot{r}_t - r_t |\dot{x}_t|^2 \right|^2 \, dt \, d \lambda_0(z_0)
    \end{equation*}
    Let us deal with each term separately so the expressions do not become unwieldly. For the first
    \begin{align*}
        &\int_\mathfrak{C} \int_0^1 \left|r_t \ddot{x}_t + 2 \dot{r}_t \dot{x}_t \right|^2 \, dt \, d \lambda_0(z_0)\\
        =& \int_\mathfrak{C} \int_0^1 r_t^2 \left| \partial_t v_t + \nabla v_t \cdot v_t + 4 \alpha_t v_t \right|^2(x_t) \, d t \, d \lambda_0(z_0)\\
        =& \int_0^1 \int_\mathfrak{C} \left| \partial_t v_t + \nabla v_t \cdot v_t + 4 \alpha_t v_t \right|^2(x_t)  \, d (r_t^2\cdot \lambda_0)(z_0) \, dt\\
        =& \int_0^1 \int \left| \partial_t v_t + \nabla v_t \cdot v_t + 4 \alpha_t v_t \right|^2(x) \, d \mu_t(x) \, dt
    \end{align*}
    We have used that $(X_t,R_t)_\# \lambda_0 = \lambda_t$ and $\mathfrak{P} \lambda_t = \mu_t$. The second term is dealt with in exactly the same way,
    \begin{align*}
        &\int_\mathfrak{C} \int_0^1 \left|\ddot{r}_t - r_t |\dot{x}_t|^2 \right|^2 \, dt \, d \lambda_0(z_0) \\
        =& \int_\mathfrak{C} \int_0^1 4r_t^2 \left| \partial_t \alpha_t + \tfrac12 \nabla \alpha_t \cdot v_t + 2 \alpha_t^2 \right|^2(x_t) \, d t \, d \lambda_0(z_0)\\
        =& \int_0^1 \int_\mathfrak{C}  4 \left| \partial_t \alpha_t + \tfrac12 \nabla \alpha_t \cdot v_t + 2 \alpha_t^2 \right|^2(x_t) \, d(r_t^2 \cdot \lambda_0)(z_0) \, d t\\
        =& \int_0^1 \int_\mathfrak{C}  4 \left| \partial_t \alpha_t + \tfrac12 \nabla \alpha_t \cdot v_t + 2 \alpha_t^2 \right|^2(x_t) \, d\mu_t(x) \, d t.
    \end{align*}

\end{proof}

\section{Cone De Casteljau's algorithm}

\thirdthm*
\begin{proof}\label{thm:cone_de_casteljau}
% Under the conditions of proposition \ref{prop:cone_decasteljau}, 
Let $\theta_{i}:=|x_{i+1}-x_i|$. Since the expressions for $w$, $u$ and $P$ are all the same with different interpolating points, we compute the derivatives of $w$ $\rho$, $r$ and $\theta$ with the suscripts removed. To compute each derivative, we just replace with the corresponding suscript.

For $\theta = |x_1-x_0|$,
\begin{align*}
\theta \dot{\theta} &= \langle x_1-x_0, \dot{x}_1-\dot{x}_0\rangle.
\end{align*}

For a point $w_0 = (1-\rho) x_0 + \rho x_1$,
\begin{align} \label{eq:position_derivative}
\dot{w}_0 &= (x_1-x_0) \dot{\rho} + (1-\rho) \dot{x}_0 + \rho \dot{x}_1.
\end{align}

For the mass $r^2 = r_0^2 (1-t)^2 + r_1^2 t^2 + 2r_0r_1t(1-t)\cos(\theta)$,
\begin{align} \label{eq:mass_derivative}
r \dot{r} =~ & 2r_0\dot{r}_0(1-t)^2 -r_0^2(1-t) + r_1\dot{r}_1t^2 + r_1^2 t \\
& + (\dot r_1 r_0 +\dot r_0 r_1)t(1-t)\cos(\theta) \notag \\
& + r_0 r_1(1-2t)\cos(\theta)  - r_0 r_1t(1-t)\sin(\theta) \dot\theta. \notag
\end{align}

For the local time $\rho = \frac1\theta \arccos\left(\frac{r_0(1-t) + r_1t\cos(\theta)}{r}\right)$,
\begin{align} \label{eq:time_derivative}
\dot\rho \theta + \rho\dot\theta =~ &  \frac{1}{r_{w_0}^2}\big( (1-t)t \sin(\theta) (\dot r_1 r_0 - \dot r_0 r_1) \\
& + (1-t)t \cos(\theta) r_0 r_1 \dot\theta \notag \\
& + r_0 r_1 \sin(\theta) + r_1^2 \dot\theta t^2 \big). \notag
\end{align}

For the first interpolation points $w_0,w_1,w_2$ we have in particular that $\dot{x}_i, \dot{r}_{x_i}= 0$, hence
\begin{align*}
\dot{\theta}_{w_i} &= 0,\\
\dot{w}_i &= (x_{i+1}-x_i) \dot{\rho}_{w_i}, \\
\dot r_{w_i} &= \frac{1}{r_{w_i}} \left(-r_{x_i}^2(1-t) + r_{x_{i+1}}^2 t  + r_{x_i} r_{x_{i+1}}(1-2t)\cos(\theta_{w_i})\right), \\
\dot\rho_{w_i} &= \frac{ r_{x_i} r_{x_{i+1}}}{\theta_{w_i}} \sin(\theta_{w_i}) - \rho_{w_i} \frac{\dot \theta_{w_i}}{\theta_{w_i}}.
\end{align*}
Thus at $t=0$,

\begin{align}
\dot{\theta}_{w_i}(0) &= 0 \\
\dot\rho_{w_i}(0) &= \frac{ r_{x_{i+1}}}{r_{x_i}} \frac{\sin(\theta_{i})}{\theta_{i}} \\
\dot r_{w_i}(0) &= r_{x_{i+1}}\cos(\theta_{i}) - r_{x_i} \\
\dot{w}_i(0) &= (x_{i+1}-x_i) \frac{ r_{x_{i+1}}}{r_{x_i}} \frac{\sin(\theta_{i})}{\theta_{i}} 
\end{align}

and $t=1$,
\begin{align}
\dot{\theta}_{w_i}(1) &= 0 \\
\dot\rho_{w_i}(1) &= \frac{r_{x_i}}{ r_{x_{i+1}}} \frac{\sin(\theta_{i})}{\theta_{i}} \\
\dot r_{w_i}(1) &= r_{x_{i+1}}-r_{x_i}\cos(\theta_{i}) \\
\dot{w}_i(1) &= (x_{i+1}-x_i) \frac{r_{x_i}}{ r_{x_{i+1}}} \frac{\sin(\theta_{i})}{\theta_{i}}.
\end{align}

We substitute these expression back into the formulas for the derivatives of position \eqref{eq:position_derivative}, mass \eqref{eq:mass_derivative} and local time \eqref{eq:time_derivative} for $u_j$ and $p$, for $j=0,1$. Notice that for both $t=0$ and $t=1$ the depdence on $\dot\theta$ vanishes.

\begin{align}
\dot\rho_{u_j}(0) &= \frac{ r_{x_{j+1}}}{r_{x_j}} \frac{\sin(\theta_{j})}{\theta_{j}} \\
\dot r_{u_j}(0) &= 2 \left( r_{x_{j+1}}\cos(\theta_{j}) - r_{j}\right) \\
\dot{u}_j(0) &= 2(x_{j+1}-x_j) \frac{ r_{x_{j+1}}}{r_{x_j}} \frac{\sin(\theta_{j})}{\theta_{j}} 
\end{align}

\begin{align}
\dot\rho_{u_j}(1) &= \frac{r_{x_j}}{ r_{x_{j+1}}} \frac{\sin(\theta_{{j+1}})}{\theta_{j+1}} \\
\dot r_{u_j}(1) &= 2 \left( r_{x_{j+1}} - r_{x_j}\cos(\theta_{j+1}) \right) \\
\dot{u}_j(1) &= 2(x_{j+1}-x_j) \frac{ r_{x_{j+1}}}{r_{x_j}} \frac{\sin(\theta_{j+1})}{\theta_{j+1}} 
\end{align}

Finally for $p$,

\begin{align}
\dot\rho_{p}(0) &= \frac{ r_{x_1}}{r_{x_0}} \frac{\sin(\theta_0)}{\theta_0}, \\
\dot r_{p}(0) &= 3 \left( r_{x_{1}}\cos(\theta_0) - r_{x_0}\right), \\
\dot{p}(0) &= 3(x_{1}-x_0) \frac{ r_{x_{1}}}{r_{x_0}} \frac{\sin(\theta_0)}{\theta_0}, 
\end{align}

\begin{align}
\dot\rho_{p}(1) &= \frac{r_{x_2}}{ r_{x_{3}}} \frac{\sin(\theta_{2})}{\theta_{2}}, \\
\dot r_{p}(1) &= 3 \left( r_{x_{3}} - r_{x_2}\cos(\theta_{2}) \right), \\
\dot{p}(1) &= 3(x_{3}-x_2) \frac{ r_{x_2}}{r_{x_3}} \frac{\sin(\theta_{2})}{\theta_{2}} .
\end{align}

\end{proof}

\bibliographystyle{apalike} 
\bibliography{references}

\begin{thebibliography}{}

\bibitem[Absil et~al., 2016]{Absil2016}
Absil, P.-A., Gousenbourger, P.-Y., Striewski, P., and Wirth, B. (2016).
\newblock Differentiable piecewise-b{\'e}zier surfaces on riemannian manifolds.
\newblock {\em SIAM Journal on Imaging Sciences}, 9(4):1788--1828.

\bibitem[Benamou and Brenier, 2000]{benamou2000computational}
Benamou, J.-D. and Brenier, Y. (2000).
\newblock A computational fluid mechanics solution to the monge-kantorovich
  mass transfer problem.
\newblock {\em Numerische Mathematik}, 84(3):375--393.

\bibitem[Benamou et~al., 2019]{benamou2019second}
Benamou, J.-D., Gallou{\"e}t, T.~O., and Vialard, F.-X. (2019).
\newblock Second-order models for optimal transport and cubic splines on the
  wasserstein space.
\newblock {\em Foundations of Computational Mathematics}, 19(5):1113--1143.

\bibitem[Chen et~al., 2018]{chen2018measure}
Chen, Y., Conforti, G., and Georgiou, T.~T. (2018).
\newblock Measure-valued spline curves: An optimal transport viewpoint.
\newblock {\em SIAM Journal on Mathematical Analysis}, 50(6):5947--5968.

\bibitem[Chewi et~al., 2020a]{chewi2020fast}
Chewi, S., Clancy, J., Gouic, T.~L., Rigollet, P., Stepaniants, G., and
  Stromme, A.~J. (2020a).
\newblock Fast and smooth interpolation on wasserstein space.
\newblock {\em arXiv preprint arXiv:2010.12101}.

\bibitem[Chewi et~al., 2020b]{chewi2020gradient}
Chewi, S., Maunu, T., Rigollet, P., and Stromme, A.~J. (2020b).
\newblock Gradient descent algorithms for bures-wasserstein barycenters.
\newblock In {\em Conference on Learning Theory}, pages 1276--1304. PMLR.

\bibitem[Chizat, 2017]{chizat2017unbalanced}
Chizat, L. (2017).
\newblock {\em Unbalanced optimal transport: Models, numerical methods,
  applications}.
\newblock PhD thesis, PSL Research University.

\bibitem[Chizat et~al., 2018]{chizat2018scaling}
Chizat, L., Peyr{\'e}, G., Schmitzer, B., and Vialard, F.-X. (2018).
\newblock Scaling algorithms for unbalanced optimal transport problems.
\newblock {\em Mathematics of Computation}, 87(314):2563--2609.

\bibitem[Gigli, 2012]{gigli2012second}
Gigli, N. (2012).
\newblock {\em Second Order Analysis on $(\mathscr{P}_2 (M), W_2)$}.
\newblock American Mathematical Soc.

\bibitem[Gousenbourger et~al., 2019]{Absil2019}
Gousenbourger, P.-Y., Massart, E., and Absil, P.-A. (2019).
\newblock Data fitting on manifolds with composite b{\'e}zier-like curves and
  blended cubic splines.
\newblock {\em Journal of Mathematical Imaging and Vision}, 61(5):645--671.

\bibitem[Kondratyev et~al., 2016]{kondratyev2016new}
Kondratyev, S., Monsaingeon, L., Vorotnikov, D., et~al. (2016).
\newblock A new optimal transport distance on the space of finite radon
  measures.
\newblock {\em Advances in Differential Equations}, 21(11/12):1117--1164.

\bibitem[Lavenant et~al., 2021]{lavenant2021towards}
Lavenant, H., Zhang, S., Kim, Y.-H., and Schiebinger, G. (2021).
\newblock Towards a mathematical theory of trajectory inference.
\newblock {\em arXiv preprint arXiv:2102.09204}.

\bibitem[Liero et~al., 2016]{liero2016optimal}
Liero, M., Mielke, A., and Savar{\'e}, G. (2016).
\newblock Optimal transport in competition with reaction: The
  hellinger--kantorovich distance and geodesic curves.
\newblock {\em SIAM Journal on Mathematical Analysis}, 48(4):2869--2911.

\bibitem[Liero et~al., 2018]{liero2018optimal}
Liero, M., Mielke, A., and Savar{\'e}, G. (2018).
\newblock Optimal entropy-transport problems and a new hellinger--kantorovich
  distance between positive measures.
\newblock {\em Inventiones mathematicae}, 211(3):969--1117.

\bibitem[Lisini, 2007]{lisini2007characterization}
Lisini, S. (2007).
\newblock Characterization of absolutely continuous curves in wasserstein
  spaces.
\newblock {\em Calculus of variations and partial differential equations},
  28(1):85--120.

\bibitem[Maniglia, 2007]{maniglia2007probabilistic}
Maniglia, S. (2007).
\newblock Probabilistic representation and uniqueness results for
  measure-valued solutions of transport equations.
\newblock {\em Journal de math{\'e}matiques pures et appliqu{\'e}es},
  87(6):601--626.

\bibitem[Otto, 2001]{otto2001geometry}
Otto, F. (2001).
\newblock The geometry of dissipative evolution equations: the porous medium
  equation.

\bibitem[Peyr{\'e} et~al., 2019]{peyre2019computational}
Peyr{\'e}, G., Cuturi, M., et~al. (2019).
\newblock Computational optimal transport: With applications to data science.
\newblock {\em Foundations and Trends{\textregistered} in Machine Learning},
  11(5-6):355--607.

\bibitem[Pooladian and Niles-Weed, 2021]{pooladian2021entropic}
Pooladian, A.-A. and Niles-Weed, J. (2021).
\newblock Entropic estimation of optimal transport maps.
\newblock {\em arXiv preprint arXiv:2109.12004}.

\bibitem[Schiebinger et~al., 2019]{schiebinger2019optimal}
Schiebinger, G., Shu, J., Tabaka, M., Cleary, B., Subramanian, V., Solomon, A.,
  Gould, J., Liu, S., Lin, S., Berube, P., et~al. (2019).
\newblock Optimal-transport analysis of single-cell gene expression identifies
  developmental trajectories in reprogramming.
\newblock {\em Cell}, 176(4):928--943.

\end{thebibliography}

\end{document}